\documentclass[11pt,a4paper]{article}
\usepackage{amsmath}
\usepackage{amsthm}
\usepackage{amsfonts}
\usepackage{amssymb}
\usepackage{stmaryrd}
\usepackage{latexsym}

\theoremstyle{definition}

\begin{document}

 \tolerance2500

\title{\Large{\textbf{On some groupoids of small orders with Bol-Moufang type of identities}}}
\author{\normalsize {Vladimir Chernov, Alexander  Moldovyan, Victor Shcherbacov}
}

 \maketitle

\begin{abstract}
We count number of groupoids of order 3  with some Bol-Moufang type identities.

\medskip


\medskip

\end{abstract}

\bigskip

\section{Introduction}

A binary groupoid $(G, \cdot)$ is  a non-empty set $G$ together with a binary operation \lq\lq $\cdot$\rq\rq. This definition is very general, therefore usually groupoids with some identities are studied. For example, groupoids with identity associativity (semi-groups) are researched.

We continue the study of groupoids with some Bol-Moufang type identities \cite{NOVIKOV_08, VD, 2017_Scerb}.
Here we present results published in \cite{CHErnov, CHErnov_2018}.

\medskip

{\bf Definition.} \label{Bol_Moufang_TYpe_Id}
Identities that involve three variables, two of which appear once on both sides of the equation and one of which appears twice on both sides are called Bol-Moufang type identities. \index{identity!Bol-Moufang type}

\medskip

Various properties of Bol-Moufang type identities in quasigroups and loops are studied in \cite{Fenyves_1, Ph_2005, Cote, AKHTAR}.

Groupoid $(Q, \ast)$ is called a quasigroup, if the following conditions are true \cite{VD}:
$(\forall u, v \in  Q) (\exists ! \, x, y \in  Q) (u * x = v  \, \& \,  y * u = v)$.

For groupoids  the following natural problems are researched: how many groupoids with some identities  of small order there exist?
A list of numbers of  semigroups of orders up to 8 is given in  \cite{Satoh}; a list of numbers of quasigroups up to 11 is given in \cite{HOP, WIKI_44}.

\section{Some results}

Original algorithm is elaborated  and corresponding program is written for generating of groupoids of small (2 and 3) orders with some Bol-Moufang identities, which are well known in quasigroup theory.

To verify the correctness of the written program the number of semigroups of order 3 was counted. Obtained result coincided with well known, namely, there exist 113 semigroups of order 3.

The following identities have the property that any of them define a commutative Moufang loop \cite{BRUCK_46, VD, HOP, 2017_Scerb} in the  class of loops: left (right) semimedial identity, Cote identity and its dual identity, Manin identity and its dual identity or in  the class of quasigroups (identity  (\ref{Comm_Muf_quas_Id}) and its dual identity).

\subsection{Groupoids with left semi-medial identity}

Left semi-medial identity in a groupoid $(Q, \ast)$ has the following form: $xx*yz=xy*xz$.
Bruck \cite{BRUCK_46, VD, 2017_Scerb} uses namely this identity to define commutative Moufang loops in the class of loops.

There exist 10 left semi-medial groupoids of order 2. There exist 7 non-isomorphic left semi-medial groupoids of order 2. The first five of them are semigroups \cite{WIKI_44}.

\[
\begin{array}{lcrr}
\begin{array}{l|ll}
\ast&1&2\\
\hline
1&1&1\\
2&1&1\\
\end{array}
&
\begin{array}{l|ll}
\star&1&2\\
\hline
1&1&1\\
2&1&2\\
\end{array}
&
\begin{array}{l|ll}
\circ&1&2\\
\hline
1&1&1\\
2&2&2\\
\end{array}
&
\begin{array}{l|ll}
\cdot&1&2\\
\hline
1&1&2\\
2&1&2\\
\end{array}
\end{array}
\]

\[
\begin{array}{lcr}
\begin{array}{l|ll}
\diamond&1&2\\
\hline
1&1&2\\
2&2&1\\
\end{array}
&
\begin{array}{l|ll}
\odot&1&2\\
\hline
1&2&1\\
2&2&1\\
\end{array}
&
\begin{array}{l|ll}
\bullet&1&2\\
\hline
1&2&2\\
2&1&1\\
\end{array}
\end{array}
\]

There exist  399 left semi-medial groupoids  of order 3.

The similar results are true for groupoids with right  semi-medial identity  $xy*zz=xz*yz$. It is clear that the identities  of left and right semi-mediality are dual. In other language they are (12)-parastrophes of each other \cite{VD, 2017_Scerb}.

It is clear that groupoids with dual identities have similar properties, including the number of groupoids of a fixed order.

\subsection{Groupoids with  Cote identity}

Identity  $x(xy*z) = (z*xx)y$ is discovered in  \cite{Cote}. Here we name this identity  Cote identity.

There exist 6 groupoids of order 2 with Cote identity. There exist 3  non-isomorphic in pairs groupoids of order 2 with Cote identity.

There exist  99 groupoids of order  3 with  Cote identity.

The similar results are true for groupoids with the following  identity  $(z\ast yx)x = y(xx \ast z)$. The last identity is (12)-parastrophe  of Cote identity.

\subsection{Groupoids with  Manin identity}

The identity $x(y*xz) = (xx*y)z$ we call  Manin identity \cite{MANIN}. The following identity is dual identity to Manin identity: $(zx\ast y)x = z(y\ast xx)$.

There exist 10 groupoids of order 2 with Manin identity. There exist 7  non-isomorphic in pairs groupoids of order 2 with Manin identity.

There exist  167 groupoids of order  3 with  Manin identity.

\subsection{Groupoids with  identity $(xy\ast x)z  = (y\ast xz) x$   \label{Comm_Muf_quas_Id} (identity (\ref{Comm_Muf_quas_Id}))}

Some properties of identity (\ref{Comm_Muf_quas_Id}) are given  in  \cite{VS_2014_Kiev, 2017_Scerb}.
 The following identity is dual identity to identity (\ref{Comm_Muf_quas_Id}): $z(x\ast yx) = x(zx\ast y)$.

There exist 6 groupoids of order 2 with identity (\ref{Comm_Muf_quas_Id}).
 There exist 3  non-isomorphic in pairs groupoids of order 2 with  (\ref{Comm_Muf_quas_Id}) identity. Any of these groupoids is a semigroup.

There exist  117 groupoids of order  3 with  identity (\ref{Comm_Muf_quas_Id}).

\subsection{Number of groupoids  of order 3 with some identities}

We count number of groupoids  of order 3 with some identities.
We use list of Bol-Moufang type identities given in \cite{Cote}. In Table 1  we present number of groupoids of order 3 with the respective identity.

\begin{table}
\centering
\caption{Number of groupoids  of order 3 with some identities.}
\footnotesize{
\[
\begin{array}{|c||c| c| c| c|}
\hline
 Name & Abbreviation & Identity  & Number  \\
\hline\hline
Semigroups &  SGR &  x(yz) = (xy)z & 113\\
\hline
 Extra &  EL &  x(y(zx)) = ((xy)z)x & 239\\
\hline
Moufang &  ML & (xy)(zx) = (x(yz))x & 196\\
\hline
Left Bol &  LB & x(y(xz)) = (x(yx))z & 215\\
\hline
Right  Bol &  RB & y((xz)x) = ((yx)z)x & 215\\
\hline
C-loops &  CL & y(x(xz)) = ((yx)x)z & 133\\
\hline
LC-loops &  LC & (xx)(yz) = (x(xy))z & 220\\
\hline
RC-loops &  RC & y((zx)x) = (yz)(xx) & 220\\
\hline
Middle Nuclear Square &  MN & y((xx)z) = (y(xx))z & 350\\
\hline
Right Nuclear Square &  RN & y(z(xx)) = (yz)(xx) & 932\\
\hline
Left Nuclear Square &  LN & ((xx)y)z = (xx)(yz) & 932\\
\hline
Comm. Moufang &  CM & (xy)(xz) = (xx)(zy) & 297\\
\hline
Abelian Group &  AG & x(yz) = (yx)z & 91\\
\hline
Comm. C-loop  &  CC & (y(xy))z = x(y(yz)) & 169\\
\hline
Comm. Alternative  &  CA & ((xx)y)z = z(x(yx)) & 110\\
\hline
Comm. Nuclear square  &  CN & ((xx)y)z = (xx)(zy) & 472\\
\hline
Comm. loops  &  CP & ((yx)x)z = z(x(yx)) & 744\\
\hline
Cheban \, 1  &  C1 & x((xy)z) = (yx)(xz) & 219\\
\hline
Cheban \, 2  &  C2 & x((xy)z) = (y(zx))x & 153\\
\hline
Lonely \, I  &  L1 & (x(xy))z = y((zx)x) & 117\\
\hline
Cheban\, I\, Dual  &  CD & (yx)(xz) = (y(zx))x  & 219\\
\hline
Lonely \, II   &  L2 & (x(xy))z = y((xx)z)  & 157\\
\hline
Lonely \, III   &  L3 & (y(xx))z = y((zx)x) & 157\\
\hline
Mate \, I   &  M1 & (x(xy))z = ((yz)x)x & 111\\
\hline
Mate \, II   &  M2 & (y(xx))z = ((yz)x)x & 196\\
\hline
Mate \, III   &  M3 & x(x(yz)) = y((zx)x) & 111\\
\hline
Mate \, IV   &  M4 & x(x(yz)) = y((xx)z) & 196\\
\hline
Triad \, I   &  T1 & (xx)(yz) = y(z(xx)) & 162\\
\hline
Triad \, II   &  T2 & ((xx)y)z = y(z(xx)) & 180\\
\hline
Triad \, III   &  T3 & ((xx)y)z = (yz)(xx) & 162\\
\hline
Triad \, IV   &  T4 & ((xx)y)z = ((yz)x)x & 132\\
\hline
Triad \, V   &  T5 & x(x(yz)) = y(z(xx)) & 132\\
\hline
Triad \, VI   &  T6 & (xx)(yz) = (yz)(xx) & 1419\\
\hline
Triad \, VII   &  T7 & ((xx)y)z = ((yx)x)z  & 428\\
\hline
Triad \, VIII   &  T8 & (xx)(yz) = y((zx)x)  & 120\\
\hline
Triad \, IX   &  T9 & (x(xy))z = y(z(xx)) & 102\\
\hline
Frute   &  FR & (x(xy))z = (y(zx))x & 129\\
\hline
Crazy Loop   & CR & (x(xy))z = (yx)(xz) & 136\\
\hline
Krypton   & KL & ((xx)y)z = (x(yz))x & 268\\
\hline
\end{array}
\]}
\end{table}

\bigskip

\textbf{Acknowledgments.} Authors thank Dr. V.D. Derech for his information on semigroups of small orders.


\begin{thebibliography}{10}

\bibitem{AKHTAR}
Reza Akhtar, Ashley Arp, Michael Kaminski, Jasmine~Van Exel, Davian Vernon, and
  Cory Washington.
\newblock \protect{The varieties of Bol-Moufang quasigroups defined by a single
  operation}.
\newblock {\em Quasigroups Related Systems}, 20(1):1--10, 2012.

\bibitem{VD}
V.D. Belousov.
\newblock {\em Foundations of the Theory of Quasigroups and Loops}.
\newblock Nauka, Moscow, 1967.
\newblock (in Russian).

\bibitem{BRUCK_46}
R.H. Bruck.
\newblock Contribution to the theory of loops.
\newblock {\em Trans. Amer. Math. Soc.}, 60:245--354, 1946.

\bibitem{CHErnov_2018}
Vladimir Chernov, Alexander Moldovyan, and Victor Shcherbacov.
\newblock \protect{On some groupoids of order three with Bol-Moufang type of
  identities}.
\newblock In {\em Proceedings of the Conference on Mathematical Foundations of
  Informatics MFOI2018, July 2-6, 2018, Chisinau}, pages 17--20, Chisinau,
  Moldova, 2018.

\bibitem{CHErnov}
Vladimir Chernov, Nicolai Moldovyan, and Victor Shcherbacov.
\newblock On some groupoids of small order.
\newblock In {\em The Fourth Conference of Mathematical Society of the Republic
  of Moldova dedicated to the centenary of Vladimir Andrunachievici
  (1917-1997), June 28 - July 2, 2017, Chisinau, Proceedings CMSM4'}, pages
  51--54, Chisinau, Moldova, 2017.

\bibitem{Cote}
B.~Cote, B.~Harvill, M.~Huhn, and A.~Kirchman.
\newblock \protect{Classification of loops of generalized Bol-Moufang type}.
\newblock {\em Quasigroups Related Systems}, 19(2):193--206, 2011.

\bibitem{Fenyves_1}
F.~Fenyves.
\newblock \protect{Extra loops. II. On loops with identities of Bol-Moufang
  type}.
\newblock {\em Publ. Math. Debrecen}, 16:187--192, 1969.

\bibitem{MANIN}
Yu.I. Manin.
\newblock {\em Cubic forms}.
\newblock Nauka, Moscow, 1972.
\newblock (in Russian).

\bibitem{NOVIKOV_08}
B.~V. Novikov.
\newblock On decomposition of moufang groupoids.
\newblock {\em Quasigroups and related systems}, 16(1):97--101, 2008.

\bibitem{HOP}
H.O. Pflugfelder.
\newblock {\em Quasigroups and Loops: Introduction}.
\newblock Heldermann Verlag, Berlin, 1990.

\bibitem{Ph_2005}
J.~D. Phillips and Petr Vojtechovsky.
\newblock \protect{The varieties of loops of Bol-Moufang type}.
\newblock {\em Algebra Universalis}, 54(3):259--271, 2005.

\bibitem{Satoh}
S.~Satoh, K.~Yama, and M.~Tokizawa.
\newblock Semigroups of order 8.
\newblock {\em Semigroup forum}, 49:7--29, 1994.

\bibitem{VS_2014_Kiev}
V.A. Shcherbacov.
\newblock \protect{About commutative Moufang loops}.
\newblock In {\em International algebraic conference dedicated to $100^{th}$
  anniversary of L.A. Kaluzhnin, July 7-12, 2014, Kyiv, Ukraine, Book of
  Abstracts}, pages 78--79, Kyiv, 2014. Taras Shevchenko State University.

\bibitem{2017_Scerb}
Victor Shcherbacov.
\newblock {\em Elements of Quasigroup Theory and Applications}.
\newblock CRC Press, Boca Raton, 2017.

\bibitem{WIKI_44}
Wikipedia.
\newblock \protect{Semigroup with two elements}, 2015.

\newblock https://en.wikipedia.org/wiki/Semigroup\_with\_two\_elements.

\end{thebibliography}

\vspace{2mm}
\begin{center}
\begin{parbox}{118mm}{\footnotesize
Vladimir Chernov$^{1}$, Nicolai  Moldovyan$^{2}$,
Victor Shcherbacov$^{3}$
\vspace{3mm}

\noindent $^{1}$Master/Shevchenko Transnistria State University

\noindent Email: volodea.black@gmail.com

\vspace{3mm}

\noindent $^{2}$Professor/St. Petersburg Institute for Informatics and Automation of Russian Academy of Sciences

\noindent Email: nmold@mail.ru

\vspace{3mm}

\noindent $^{3}$Principal Researcher/Institute of Mathematics and Computer Science of Moldova

\noindent Email: victor.scerbacov@math.md

}

\end{parbox}
\end{center}

\end{document}